\documentclass[12pt,twoside,centertags,reqno]{article}
\usepackage{amsfonts,amssymb,amsmath,amsthm,latexsym,amsxtra,amscd}
\usepackage[english]{babel}
\usepackage{verbatim,alltt}
\usepackage{epsfig,epic,eepic}
\usepackage{fancyheadings}
\usepackage{xypic}
\newtheorem{lemma}{Lemma}[section]
\newtheorem{proposition}[lemma]{Proposition}
\newtheorem{theorem}[lemma]{Theorem}

\newtheorem{conjecture}[lemma]{Conjecture}
\theoremstyle{definition}
\newtheorem{definition}[lemma]{Definition}

\date{}
\newcommand{\addresses}[3]{}
\newcommand{\emails}[3]{}
\newcommand{\classification}[1]{%
\renewcommand{\thefootnote}{}%
\footnotetext{\mbox{\hspace*{-12pt}
1991 {\it Mathematics Subject Classification.}}
{#1}%
}}
\newcommand{\keywords}[1]{\renewcommand{\thefootnote}{}%
\footnotetext{\mbox{\hspace*{-12pt}
{\it Key words.}} 
{#1}%
}}
\title{{\bf Monodromy}}
\author{Wolfgang Ebeling}
\addresses{
Universit\"{a}t Hannover, Institut f\"{u}r Mathematik,
Postfach 6009, D-30060 Hannover, Germany 
}
\emails{ebeling@math.uni-hannover.de
}
\newcommand{\CC}{{\Bbb C}}

\newcommand{\QQ}{{\Bbb Q}}

\newcommand{\RR}{{\Bbb R}}
\newcommand{\ZZ}{{\Bbb Z}}

\newcommand{\NN}{{\Bbb N}}

\newcommand{\eps}{\varepsilon}

\begin{document}
\maketitle
\begin{center}
{\it Dedicated to Gert-Martin Greuel on the occasion of his 60th birthday.}
\end{center}
\thispagestyle{empty}
\classification{14D05, 32S40
%
}
\keywords{monodromy, zeta function, spectrum, isolated singularity
%
}
\renewcommand{\thefootnote}{\arabic{footnote}}%

\begin{abstract}
\noindent
Let $(X,x)$ be an isolated complete intersection singularity and let $f : (X,x)
\to (\CC,0)$ be the germ of an analytic function with an isolated singularity
at $x$. An important topological invariant in this situation is the
Picard-Lefschetz monodromy operator associated to $f$. We give a survey on what
is known about this operator. In particular, we review methods of computation of the
monodromy and its eigenvalues (zeta function), results on the Jordan normal
form of it, definition and properties of the spectrum, and the relation between the monodromy and
the topology of the singularity.
%
%
\end{abstract}

%
%



\section*{Introduction}
The word 'monodromy' comes from the greek word
${\mu}o{\nu}o-{\delta}{\rho}o{\mu}{\psi}$ and means something like
'uniformly running' or 'uniquely running'. According to
\cite[3.4.4]{Lamotke05}, it was first used by B.~Riemann \cite{Riemann57}.
It arose in keeping track of the solutions of the hypergeometric
differential equation going once around a singular point on a closed path
(cf.\ \cite{Brieskorn76}). The group of linear substitutions which the
solutions are subject to after this process is called the {\em monodromy
group}.

Since then, monodromy groups have played a substantial r\^ole in many areas
of mathematics. As is indicated on the webside 'www.monodromy.com' of
N.~M.~Katz, there are several incarnations, classical and $l$-adic, local and
global, arithmetic and geometric. Here we concentrate on the classical local
geometric monodromy in singularity theory. More precisely we focus on the
monodromy operator of an isolated hypersurface or complete intersection
singularity. The investigation of this operator started in 1967 with the
proof of the famous monodromy theorem (see \S
1). This theorem can be proved using the theory of the Gau\ss-Manin connection which was introduced
by E.~Brieskorn for isolated hypersurface singularities \cite{Brieskorn70}. The study of this
connection for isolated complete intersection singularities was started by
G.-M.~Greuel in his thesis \cite{Greuel75}. 

We try to review
the results of 37 years of investigation of the monodromy operator.  The results
include results on the zeta function of the monodromy and the spectrum of a
singularity. The monodromy contains a lot of information
about the topology of the singularity. This was one motivation to study the monodromy. We review the
known facts in the last section.

For a basic introduction to the subject for non-specialists see
\cite{Ebeling01}.

Aspects which are not mentioned or only touched in this survey are
\begin{itemize}
\item Monodromy of a polynomial function. For a survey on this topic
see \cite{Dimca00} and \cite{GLM01}.
\item Generalizations to non-isolated singularities. Here we refer to the
survey  of D.~Siersma \cite{Siersma01}.
\item Monodromy groups. We do not talk about monodromy groups of isolated complete intersection
singularities. We refer to our book \cite{Ebeling87} for this topic.
\item Braid monodromy. 
Recently, M.~L\"onne introduced a notion of braid monodromy of singularities.  
He computed the braid monodromy and the fundamental group
of the complement of the discriminant of a Brieskorn-Pham singularity
\cite{Lonne03},  making a
substantial contribution to the last problem \cite[Probl\`eme
20]{Brieskorn73} of Brieskorn's list of problems on monodromy.
\end{itemize}

\section{The monodromy operator}
Let $(Y,0) \subset (\CC^N,0)$ be an isolated complete intersection
singularity (abbreviated ICIS in the sequel) of dimension $n+1$, i.e. $(Y,0)$ is the germ of an
analytic variety of pure dimension $n+1$ with an isolated singularity at the origin given by
$Y=F^{-1}(0)$, where $F=(f_1, \ldots , f_{N-n-1})
: (\CC^N,0) \to (\CC^{N-n-1},0)$ is the germ of an analytic mapping. Let $f: \CC^N \to
\CC$ be an analytic function such that the restriction $f: Y \to \CC$ which we denote by the
same symbol has an isolated
singularity at the origin. We assume that
$f(0)=0$. Let $\eps
>0$ be small enough such that the closed ball $B_\eps \subset \CC^N$ of
radius $\eps$ around the origin in $\CC^N$ intersects the fibre $f^{-1}(0)$
transversely. Let $0 < \delta \ll \eps$ be such that for $t$ in the disc
$D_\delta \subset \CC$ around the origin, the fibre $f^{-1}(t) \cap Y$
intersects the ball $B_\eps$ transversely. Let 
\begin{eqnarray*} 
X_t & := & f^{-1}(t) \cap B_\eps \cap Y \mbox{ for } t \in D_\delta,\\ 
X & := & f^{-1}(D_\delta) \cap B_\eps \cap Y, \\ 
X'& := & X \setminus X_0,\\ 
D' & := & D \setminus \{ 0\}.
\end{eqnarray*} 
Then $(X_0,0)$ is an ICIS of dimension
$n$. In the important special case when $Y$ is smooth,
$(Y,0)=(\CC^{n+1},0)$, then $(X_0,0)$ is an isolated hypersurface
singularity. By a result of J.~Milnor \cite{Milnor68} in the case when $Y$
is smooth and H.~Hamm \cite{Hamm71} in the general case, the mapping
$f|_{X'} : X' \to D'$ is the projection of a locally trivial $C^\infty$
fibre bundle. A fibre $X_t$ of this bundle is called {\em Milnor fibre}. It
has the homotopy type of a bouquet of $\mu$ $n$-spheres where $\mu$ is the
Milnor number. Therefore its only interesting homology group is the group
$H_n(X_t, \ZZ)$. It is of rank $\mu$. Parallel translation along the path
$$\gamma: [0,1] \to D_\delta, \quad t \mapsto \delta e^{2 \pi i t},$$
yields a diffeomorphism $h : X_\delta \to X_\delta$ called the {\em
geometric monodromy} of the singularity. It is determined up to isotopy.

\begin{definition} The induced homomorphism $h^\CC_\ast : H_n(X_\delta,\CC) \to
H_n(X_\delta,\CC)$ (resp.\ $h^\ZZ_\ast : H_n(X_\delta,\ZZ) \to H_n(X_\delta,\ZZ)$) is
called the {\em complex} (resp.\ {\em integral}) {\em monodromy (operator)}
of the singularity.
\end{definition}

This operator is also sometimes called the {\em Picard-Lefschetz monodromy
operator} since the
consideration of this operator goes back to E.~Picard \cite{PS97} and
S.~Lefschetz
\cite{Lefschetz24} (see also \cite{Lamotke75}, \cite{Pham65}).

\begin{theorem}[Monodromy theorem] \label{Thm1}
\begin{itemize}
\item[{\rm (a)}] The eigenvalues of $h_\ast$ are roots of unity.
\item[{\rm (b)}] The size of the blocks in the Jordan normal form of
$h_\ast$ is at most $(n+1) \times (n+1)$.
\item[{\rm (c)}] If $(Y,0)$ is smooth, then the size of the Jordan blocks
for the eigenvalue $1$ is at most $n \times n$.
\end{itemize}
\end{theorem}

\begin{sloppypar}

There are many different proofs of
this theorem: by A.~Borel (unpublished), E.~Brieskorn \cite{Brieskorn70}
(for $(Y,0)$ smooth, generalized by G.-M.~Greuel \cite{Greuel75}),
C.~H.~Clemens
\cite{Clemens69}, P.~Deligne \cite{Deligne70, DK73}, P.~A.~Griffiths
\cite{Griffiths70},
A.~Grothendieck \cite{Grothendieck72},
N.~M.~Katz \cite{Katz70}, A.~Landman \cite{Landman67, Landman73}, L\^e D\~ung Tr\'ang \cite{Le78}
(of (a), for (b) see the book of E.~Looijenga \cite{Looijenga84}),
B.~Malgrange \cite{Malgrange74b}, and W.~Schmid \cite{Schmid73} (see also
the survey
\cite{Griffiths70}). 

Examples of B.~Malgrange \cite{Malgrange73} show that the bounds on the
size of the Jordan blocks are sharp.

For weighted homogeneous singularities with $(Y,0)$ smooth, Milnor
\cite{Milnor68} has shown that the complex monodromy $h^\CC_\ast$ is
diagonalizable.  For
weighted homogeneous ICIS this was shown by A.~Dimca \cite{Dimca85}.
For irreducible plane curve singularities, L\^e \cite{Le72} has shown
that the monodromy is of finite order. N.~A'Campo \cite{A'Campo73b} has
shown that for isolated plane
curve singularities with more than one branch the monodromy is in general
not of finite order.
A.~H.~Durfee \cite{Durfee75b} has given a necessary and sufficient
condition for the monodromy of a
degenerating family of curves to be of finite order.

We now mention several results which are only valid in the case
when $(Y,0)=(\CC^{n+1},0)$. J.~Scherk \cite{Scherk80} has shown that if
 $f^{r+1}$ belongs to the ideal $( \partial f/\partial x_0,
\ldots , \partial f/\partial x_n )$ of the ring ${\cal O}_{n+1}$ of germs
of holomorphic functions on $\CC^{n+1}$, then the size of the Jordan blocks
of $h^\CC_\ast$ is at most $(r+1) \times (r+1)$. By a theorem of
J.~Brian\c{c}on
and H.~Skoda \cite{BS74}, $f^{n+1} \in ( \partial f/\partial x_0,
\ldots , \partial f/\partial x_n )$. Therefore Scherk's theorem implies the
Monodromy theorem. Generalizations of Scherk's theorem can be found in
\cite{SchSt85}.

\end{sloppypar}

M.~G.~M.~van Doorn and J.~H.~M.~Steenbrink \cite{vDSt89} have proved the
following supplement to the Monodromy theorem: If there exists a Jordan
block of size $(n+1) \times (n+1)$, then there exists a Jordan block of
size $n \times n$ for the eigenvalue 1. Since a plane curve singularity is
reducible if and only if $h^\CC_\ast$ has an eigenvalue 1, this implies
L\^e's theorem.

Let $f: (\CC^{n+1},0) \to (\CC,0)$ and $g :(\CC^{m+1},0) \to
(\CC,0)$ be two germs of analytic functions with an isolated
singularity at 0. Denote by $c_f$ and $c_g$ the complex
monodromy operators of $f$ and $g$ respectively. Denote by
$c_{f+g}$ the complex monodromy operator of the germ $f+g$. The famous
theorem of M.~S\'ebastiani and R.~Thom \cite{ST71} states that
$$c_{f+g} = c_f \otimes c_g.$$
The author and Steenbrink \cite{ESt98} have proved a generalization of this
theorem for a suspension of an ICIS.

In the case when $Y=\CC^{n+1}$ one can associate to $f \in \CC\{x_0, \ldots
, x_n\}$ its {\em
Bernstein-Sato polynomial}. This is defined as follows. Let $s$ be a new
variable. Then there
exists a differential operator $P=P(x,s, \partial /\partial x)$ whose
coefficients are convergent
power series in $s$ and $x_0$, \dots , $x_n$ and a nonzero polynomial $b(s)
\in \CC[s]$ satisfying
the formal identity
$$P f^s = b(s) f^{s-1}.$$
The set of all polynomials $b(s) \in \CC[s]$ for which such an identity
holds (for some operator
$P$) forms an ideal, and the unique monic generator for this ideal is
called the {\em Bernstein-Sato
polynomial} of $f$. It is denoted by $b_f(s)$. According to Malgrange
\cite{Malgrange74a,
Malgrange75} (see also \cite{Bony76}) there is the following relation to
the monodromy of $f$: The
zeros $s_1$, $s_2$, \dots of
$\widetilde{b}_f(s):=b_f(s)/s$ are rational and less than 1, the minimal
polynomial of the monodromy
divides the polynomial $p(t):=\prod (t- \exp(-2 \pi i s_j))$, and on the
other hand, $p(t)$ divides
the characteristic polynomial of the monodromy. D.~Barlet \cite{Barlet86}
has shown that if there
exists a $k \times k$ Jordan block for the eigenvalue $\exp(-2 \pi i u)$ of
the monodromy, then
there exist at least $k$ (counted with multiplicity) roots of $b_f(s)$ of
the form $-q-u$ with $q
\in [0,n]$ ($0 \leq u <1$).

\section{Computation of monodromy}
We now review methods to compute the monodromy operator.

We first consider the complex monodromy operator. In the case when $(Y,0)$
is smooth, Brieskorn \cite{Brieskorn70} has indicated a method to compute
the complex monodromy operator. This method has been implemented by
M.~Schulze to {\sc Singular} \cite{Schulze99} (see also \cite{Schulze03}).

For plane curve singularities there is an algorithm in the book of
D.~Eisenbud and W.~Neumann \cite{EN85} to compute the Jordan normal form of
the complex monodromy operator from a splicing diagram of the singularity.

For superisolated surface singularities (see the article of E.~Artal-Bartolo, I.~Luengo and
A.~Melle-Hern\'andez in this volume \cite{ALM05}), Artal-Bartolo \cite{ArtalBartolo91a,
ArtalBartolo94} has determined the Jordan normal form of the complex
monodromy operator.

We now consider the integral monodromy operator. Let $Y_\eta:=F^{-1}(\eta) \cap
B_\eps$ where $\eta$ is a regular value of $F$ sufficiently close to 0.
Let $c:=h^\ZZ_\ast$ be the integral monodromy operator. The above path
$\gamma$ also induces a map $\widehat{c}: H_{n+1}(Y_\eta,X_\delta) \to
H_{n+1}(Y_\eta,X_\delta)$ on the relative homology groups. (Unless otherwise
stated, we consider homology with integral coefficients). Then we have the
following diagram with exact rows and commutative squares (cf.\ \cite{Ebeling87}):
$$\diagram
0 \rto  & H_{n+1}(Y_\eta) \rto \dto_{{\rm id}} & H_{n+1}(Y_\eta,X_\delta) \rto
\dto_{\widehat{c}} & H_n(X_\delta) \rto \dto_c & 0 \\
0 \rto  & H_{n+1}(Y_\eta) \rto & H_{n+1}(Y_\eta,X_\delta) \rto & H_n(X_\delta) \rto  & 0
\enddiagram
$$
Let $A$ be the intersection matrix on $H_{n+1}(Y_\eta,X_\delta)$ with respect to a
distinguished basis of thimbles (cf.\ \cite{Ebeling87}). It is a $\nu \times \nu$-matrix where
$\nu$ is the number of thimbles in a distinguished basis of thimbles. One has
$\nu = \mu + \mu'$ where $\mu'$ is the Milnor number of the singularity $(Y,0)$. The matrix $A$ is
encoded in the {\em Coxeter-Dynkin diagram} of the singularity.
This matrix is of
the form $A=V+(-1)^n V^t$ for some upper triangular matrix 
$$V = \left(
\begin{array}{ccccc} (-1)^{\frac{n(n+1)}{2}} & \ast & \cdots & \ast & \ast
\\ 0 & (-1)^{\frac{n(n+1)}{2}} & \ddots & \vdots & \vdots \\ \vdots &
\ddots & \ddots & \ddots & \vdots \\ 0 & \cdots & \ddots &
(-1)^{\frac{n(n+1)}{2}} & \ast \\ 0 & \cdots & \cdots & 0 &
(-1)^{\frac{n(n+1)}{2}}  \end{array} \right).$$
If $Y=\CC^{n+1}$ then $V$ is the matrix of the {\em (integral) Seifert
form} or of the inverse of
the {\em variation operator} of the singularity (see \cite{AGV88, Durfee74,
Lamotke75}). For general
$(Y,0)$, the author and S.~M.~Gusein-Zade defined in \cite{EG99} a
variation operator the inverse of
which has the matrix $V$.
The operator $\widehat{c}$ is the product of the Picard-Lefschetz
transformations corresponding to
the elements of a distinguished basis of thimbles (cf.\ \cite{Ebeling87}).
In the case when $Y=\CC^3$
and
$f$ defines a simple singularity, the Coxeter-Dynkin diagram is the
classical Coxeter-Dynkin diagram
of a root system of type $A_\mu$, $D_\mu$, $E_6$, $E_7$, or $E_8$ and
$\widehat{c}=c$ is the
corresponding {\em Coxeter element}.
It follows from \cite[Chap.~V, \S 6, Exercice 3]{Bourbaki68} (see also
\cite{Levine66}) that the
matrix
$\widehat{C}$ of the operator $\widehat{c}$ is given by
$$\widehat{C} = (-1)^{n+1} V^{-1} V^t.$$

F.~Lazzeri \cite{Lazzeri73} (see also \cite{Lazzeri88}) and independently
A.~M.~Gabrielov
\cite{Gabrielov74a} showed that
in the case when
$Y=\CC^{n+1}$ the Coxeter-Dynkin diagram is connected thus extending an
earlier result of C.~H.~Bey
\cite{Bey72a, Bey72b} for curves. A.~Hefez and Lazzeri
\cite{HL74} computed the intersection matrix of Brieskorn-Pham
singularities solving in this way an
open problem stated by Brieskorn \cite{Brieskorn73} and F.~Pham \cite{Pham65}.

Gabrielov
\cite{Gabrielov73, Gabrielov74b} has given methods to compute the
intersection matrix $A$ for some
special singularities. N.~A'Campo \cite{A'Campo75b, A'Campo75c} and
independently
Gusein-Zade \cite{GuseinZade74a, GuseinZade74b} have found a
beautiful method to compute the intersection matrix for isolated plane
curve singularities using real morsifications. This method was generalized
in \cite{EG96, EG98} to suspensions of fat points.

A rather general
method to compute an intersection matrix for isolated hypersurface
singularities using polar curves was found by Gabrielov \cite{Gabrielov79}.
This method was generalized to ICIS in \cite{Ebeling87}. The author
\cite{Ebeling80, Ebeling83} has
computed the characteristic polynomial of the monodromy for the uni- and
bimodal hypersurface singularities in Arnold's classification
\cite{Arnold75} and the
intersection matrix $A$ for the elliptic hypersurface singularities
\cite{Ebeling85}. Gusein-Zade \cite{GuseinZade76} gave a recursive formula
for the characteristic
polynomials of the monodromy for the singularities of Arnold's series of
singularities
\cite{Arnold75}.

In
\cite{Ebeling87} the matrices
$A$ were computed for the simple space curve singularities classified by
M.~Giusti
\cite{Giusti77} except $Z_9$ and $Z_{10}$ and many of the $\cal K$-unimodal
isolated singularities of complete intersection surfaces classified by
C.~T.~C.~Wall \cite{Wall83}. The missing cases $Z_9$ and $Z_{10}$ were
studied in \cite{EG95} and the case $I_{1,0}$ was treated in
\cite{Ebeling99}.

P.~Orlik and R.~Randell \cite{OR77} computed the integral monodromy for
some classes of
weighted homogeneous singularities.

If $Y=\CC^{n+1}$ and $f$ is a real analytic function, i.e.\ takes real
values on $\RR^{n+1} \subset
\CC^{n+1}$, then Gusein-Zade \cite{GuseinZade84} showed that the integral
monodromy is the product
of two involutions (see also \cite{A'Campo03}).


\section{Zeta function}
The {\em zeta function} of the monodromy is defined to be
$$\zeta(t) := \prod_{q \geq 0} \left\{ \det ({\rm id}_\ast - t h_\ast;
H_q(X_\delta,\CC)) \right\}^{(-1)^{q+1}}.$$
The relation with the characteristic polynomial $\Delta(t)$ of the monodromy
is
$$\Delta(t):= \det(t {\rm id}_\ast - h_\ast) = t^\mu \left[ \frac{t-1}{t}
\zeta \left( \frac{1}{t} \right) \right]^{(-1)^{n+1}}.$$
The {\em Lefschetz numbers} are defined by
$$\Lambda_k := \Lambda(h^k_\ast) = \sum_{ q \geq 0} (-1)^q {\rm
Tr}[h_\ast^k; H_q(X_\delta,\CC)].$$
We define rational numbers $\chi_m$ by
$$\Lambda_k = \sum_{m | k} m \chi_m.$$
Explicitely, these numbers can be defined by M\"obius inversion
$$\chi_m = \frac{1}{m} \left( \sum_{k|m} \mu \left( \frac{m}{k} \right)
\Lambda_k \right),$$
where $\mu( \ )$ denotes the M\"obius function.
By A.~Weil (cf.\ \cite{Milnor68}) we have
$$\zeta(t) = \prod_{m \geq 1} (1-t^m)^{-\chi_m}.$$

The following statements were explained to me by D.~Zagier.

\begin{proposition}[Zagier]   \label{Prop1}
\begin{itemize}
\item[{\rm (i)}] The numbers $\chi_m$ are integers.
\item[{\rm (ii)}] The following statements are equivalent:
\begin{itemize}
\item[{\rm (a)}] $\Delta(t)$ is a product of cyclotomic polynomials.
\item[{\rm (b)}] $\chi_m \neq 0$ for only finitely many $m$.
\item[{\rm (c)}] The sequence $(\Lambda_k)$ is periodic.
\end{itemize}
\end{itemize}
\end{proposition}

\begin{proof} (i) is proved by induction: Assume that $\chi_m \in \ZZ$ for
$m \leq \ell$. Then
$$\prod_{m=1}^\ell (1- t^m)^{-\chi_m}$$
is a formal power series with integer coefficients which starts
with 1. Since $\zeta(t)$ is also a power series with integral coefficients,
the same is true for
$$\frac{\zeta(t)}{\prod_{m=1}^\ell (1- t^m)^{-\chi_m}} =
\prod_{m=\ell+1}^\infty (1- t^m)^{-\chi_m}.$$
But this power series starts with $1+\chi_{\ell +1} t^{\ell +1}$.

The proof of (ii) is done in several steps:

The implication (a) $\Rightarrow$ (b) is clear.

(b) $\Rightarrow$ (c): Let $\chi_m =0$ for all $m$ which do not divide a
number $Q$. Then for positive integers $d,r$ with $0 < r \leq Q$ we have
$$\Lambda_{dQ+r}= \sum_{m|dQ+r} m \chi_m = \sum_{m|r} m \chi_m =
\Lambda_r.$$

(c) $\Rightarrow$ (a): Let the sequence $(\Lambda_k)$ be periodic of period
$Q$. Let $s_k := {\rm Tr}\, h^k_\ast$ and let $p(t):= \det( {\rm id}_\ast -
t h_\ast)$. Then
$$p(t) = \exp ( {\rm Tr} (\log ({\rm id}_\ast - t h_\ast))) =
 \exp \left( - \sum_{k=1}^\infty s_k \frac{t^k}{k} \right).$$
For the logarithmic derivative of $p(t)$ we get from this
$$\frac{p'(t)}{p(t)} = - \sum_{k=1}^\infty s_k t^{k-1}.$$
Since $s_k$ is periodic of period $Q$, we get
\begin{eqnarray*}
\frac{p'(t)}{p(t)} & = & - \sum_{k=1}^Q s_k t^{k-1} (1+t^Q+t^{2Q}+ \ldots)
\\
 & = & {} - \sum_{k=1}^Q s_k \frac{t^{k-1}}{1-t^Q} = \frac{q(t)}{1-t^Q}
\end{eqnarray*}
for some polynomial $q(t)$. But each zero of $p(t)$ must be a simple pole
of the logarithmic derivative and hence a zero of $1-t^Q$.
\end{proof}

If $(Y,0)$ is smooth and $f$ is singular, then by \cite{A'Campo73a}
$\Lambda_1=0$. L\^e \cite{Le75} proved that in this case there exists a
characteristic
diffeomorphism $h$ without fixed points. In a letter to A'Campo, P.~Deligne
showed that more
generally in the case when $(Y,0)$ is smooth,
$\Lambda_k=0$ for $0 < k < {\rm mult}(f)$ where ${\rm mult}(f)$ is the
multiplicity of $f$.
G.~G.~Il'yuta \cite{Il'yuta87} gave formulae expressing the Lefschetz
numbers in terms of cycles of the Coxeter-Dynkin diagram. In
\cite{Ebeling96} the author showed that ${\rm Tr}\, C^2 = (-1)^r$ where $r$
is the corank of $f$.

Let $\pi : \widetilde{Y} \to Y$ be a resolution of $Y$ and let
$\widetilde{f}:= f \circ \pi$ be the composition. We denote by
$\widetilde{X}_0$ the proper transform of $X_0=f^{-1}(0)$. Let
$\widetilde{f}^{-1}(0) = \widetilde{X}_0 \cup E_1 \cup \cdots \cup E_s$
where $E_i$ is irreducible. We assume that the following conditions are
satisfied:
\begin{itemize}
\item[(1)] $\pi : \widetilde{X}_0 \to X_0$ is a resolution of $X_0$.
\item[(2)] Each exceptional divisor $E_i$ is smooth and
$\widetilde{f}^{-1}(0)$ has only normal crossings.
\end{itemize}
Let $m_i$ be the order of the function $\widetilde{f}$ along the divisor
$E_i$ and let
$$E'_i := E_i \setminus \left( \bigcup_{j \neq i} E_j \right) \cup
\widetilde{X}_0.$$

Then we have the following famous theorem of A'Campo \cite{A'Campo75a}:

\begin{theorem}[A'Campo]
Under the above assumptions, we have
$$\zeta(t) = \prod_{i=1}^s (1- t^{m_i})^{-\chi(E'_i)}$$
where $\chi(E'_i)$ is the topological Euler characteristic of $E'_i$.
\end{theorem}

The theorem was formulated by A'Campo only for the case $Y=\CC^{n+1}$ but
it can be easily generalized to this more general situation (see e.g.\
\cite{Oka90}). From Proposition~\ref{Prop1} we see that A'Campo's theorem
implies Part (a) of the Monodromy theorem. In \cite{AGV88} it is shown that
Part (b) can also be derived from that theorem.

A generalization of A'Campo's theorem using partial resolutions was given
by Gusein-Zade, Luengo and Melle-Hern\'andez \cite{GLM97}.

If $Y=\CC^{n+1}$ and $f$ is non-degenerate with respect to the Newton
diagram, then A.~Varchenko \cite{Varchenko76} (and also independently
F.~Ehlers \cite{Ehlers78}) have given a formula to compute $\zeta(t)$ from
the Newton diagram. This formula was generalized to the general case by
M.~Oka \cite{Oka90}.

A.~Campillo, F.~Delgado and Gusein-Zade \cite{CDG99} have shown that for an
irreducible curve singularity, the zeta function $\zeta(t)$ coincides with
the Poincar\'e series $P(t)$ of the natural filtration on the ring of
functions of such a singularity given by the order with respect to a
uniformization.

Now suppose $Y=\CC^{n+1}$ and $f$ is weighted homogeneous of weights $q_0,
\ldots , q_n$ and degree $d$. Here $q_0, \ldots , q_n$ are assumed to be
coprime. Then the geometric monodromy $h: X_1 \to X_1$ can be described as
follows \cite{Milnor68}:
$$ h(z_0, \ldots , z_n) = (e^{2 \pi i/q_0} z_0,
\ldots , e^{2 \pi i/q_n}z_n).$$
The monodromy operator $h_\ast$ is of order $d$.
Milnor and P.~Orlik \cite{MO70} have shown
how to compute $\zeta(t)$ from the weights and the degree of $f$. Greuel
and H.~Hamm \cite{GH78} have given a more general formula for a weighted
homogeneous ICIS $(Y,0)$.

K.~Saito \cite{KSaito88} has shown that if $Y=\CC^3$ then all the primitive
$d$-th roots of unity are eigenvalues of $h_\ast$.

Saito \cite{Saito98a, Saito98b} also introduced a duality between rational
functions of the form of the zeta function.
If $\phi(t)$ is a rational function of the form
$$\phi(t) = \prod_{m|d} (1-t^m)^{\chi_m} \quad \mbox{for } \chi_m \in \ZZ
\mbox{ and some } d \in \NN,$$
then he defines
$$\phi^\ast(t)= \prod_{k|d} (1-t^k)^{-\chi_{d/k}}.$$
Let $f$ be a function defining one
of the 14 exceptional unimodal hypersurface singularities in the sense of
V.~I.~Arnold \cite{Arnold76}. Arnold has observed a strange duality between
these singularities \cite{Arnold75}. Saito has observed the following fact: If $\Delta(t)$ is the
characteristic polynomial of the monodromy of $f$ then $\Delta^\ast(t)$ is the
characteristic polynomial of the monodromy of the dual singularity. The
author and C.~T.~C.Wall \cite{EW85} have found an extension of Arnold's
strange duality embracing also ICIS. The author \cite{Ebeling99} has shown
that Saito's duality also holds for this extension and he has related it to
polar duality and to a duality of weight systems found by M.~Kobayashi
\cite{Ebeling00, Ebeling05}. 

If $n=2$ and $Y=\CC^3$ or $(Y,0)$ is a certain
special ICIS, then it was shown \cite{Ebeling02} that the Saito dual
$\Delta^\ast(t)$
of the characteristic polynomial of the monodromy is equal to the product
of the Poincar\'e series $P(t)$ of the coordinate algebra and some rational
function ${\rm Or}(t)$ depending only on the orbit invariants of the natural
$\CC^\ast$-action on the singularity. In the case when
$Y=\CC^{n+1}$ and $f$ is a Newton non-degenerate function, the author and
Gusein-Zade \cite{EG02} showed that the same holds for the Saito dual of
the inverse of the reduced zeta function $\widetilde{\zeta}(t)$ (reduced means
considering reduced homology). Finally, in \cite{EG04} this was generalized
to the case when $Y$ is a complete intersection given as the zero set of
functions $f_1, \ldots, f_{k-1}$ and $f=f_k$ to the product of the Saito
duals of the inverse reduced zeta functions $\widetilde{\zeta}_j(t)$ of the
monodromy operators of $f_j$ on $f_1= \ldots = f_{j-1}=0$ for $j=1, \ldots,
k$. J.~Stevens \cite{Stevens03} proved that this result implies the theorem
of Campillo, Delgado and Gusein-Zade \cite{CDG99}.

Now let $Y=\CC^{n+1}$ and assume that $f \in \ZZ[x_0, \ldots, x_n]$. For 
a prime number $p$ denote by $\ZZ_p$ the $p$-adic integers. Consider the $p$-adic integral
$$Z_p(s):= \int_{\ZZ_p^{n+1}} \vert f(x) \vert_p^s \vert dx \vert$$
for $s \in \CC$, ${\rm Re}(s) > 0$, where $\vert dx \vert$ denotes the Haar measure on $\QQ_p^{n+1}$
normalized in such a way that $\ZZ_p^{n+1}$ is of volume 1. This function 
is called the {\em $p$-adic Igusa zeta function}.  Now there is
the following famous conjecture \cite{Igusa75} (see also \cite{Denef91}):

\begin{conjecture}[Igusa's monodromy conjecture] For almost all prime numbers $p$, if $s_0$ is a pole
of $Z_p(s)$  then $e^{2 \pi i {\rm Re}(s_0)}$ is an eigenvalue of the monodromy
operator $h_\ast$ at some point of $\{f=0\}$.
\end{conjecture}

J.~Denef and F.~Loeser \cite{DL92} defined a topological zeta function
$Z_{\rm top}(t)$ generalizing Igusa's zeta function. The analogous conjecture is
stated for this function. Loeser
\cite{Loeser88, Loeser90}, W.~Veys \cite{Veys93}, Artal-Bartolo,
P.~Cassou-Nogu\`es, Luengo and Melle-Hern\'andez \cite{ACLM02, ACLM03} and
B.~Rodrigues and Veys \cite{RV03} proved various special cases of this
conjecture. See \cite{Veys03} for an
excellent survey on this topic and the article of Artal-Bartolo, Luengo and
Melle-Hern\'andez in this volume \cite{ALM05}.

In \cite{GLM04} a motivic version of the zeta function of the monodromy is
discussed and compared with the motivic zeta function of Denef and Loeser.

\section{Spectrum}
In the case $(Y,0)$ smooth, Steenbrink \cite{Steenbrink77b} showed that
there exists a mixed Hodge structure on the Milnor fibre. Let
$H=H^n(X_\delta,\ZZ)$. Such a mixed Hodge structure consists of an increasing
{\em weight} filtration
$$0=W_{-1} \subset W_0 \subset \cdots \subset W_{2n}=H \otimes \QQ$$
of $H \otimes \QQ$ and a decreasing {\em Hodge} filtration
$$H \otimes \CC = F^0 \supset F^1 \supset \cdots \supset F^n \subset
F^{n+1}=0.$$
It follows from  a result of M.~Saito \cite{MSaito90} that in the general
case, the analogue in cohomology of the short exact sequence in \S 2 can be
considered as a sequence of mixed Hodge structures (see \cite{ESt98}).

The mixed Hodge structure is used to define the spectrum of a singularity.
The spectrum was defined by Steenbrink \cite{Steenbrink77b} and
Arnold \cite{Arnold81} in the case when $(Y,0)$ is smooth and  in
\cite{ESt98} in the general case.

\begin{definition} The {\em spectrum} ${\rm Sp}(f)$ of $f$ is defined as
follows. Let $p \in \ZZ$, $0 \leq p \leq n$. A rational number $\alpha \in
\QQ$ with $n-p-1 < \alpha \leq n-p$ is in ${\rm Sp}(f)$ if and only if
$e^{2 \pi i \alpha}$ is an eigenvalue of the semisimple part of $h^\ast$ on
$F^pH/F^{p+1}H$. Here $H=H^n(X_\delta,\CC)$ if $Y=\CC^{n+1}$ and
$H=H^{n+1}(Y_\eta,X_\delta,\CC)$ in the general case. The multiplicity of $\alpha$
is the dimension of the corresponding eigenspace.
\end{definition}

The spectrum is an unordered tuple of $\nu$ rational numbers $\alpha_1$,
\dots , $\alpha_\nu$ which lie between $-1$ and $n$. We order these numbers
as follows:
$$-1 < \alpha_1 \leq \alpha_2 \leq \ldots \leq \alpha_\nu < n.$$
There is a symmetry property
$$\alpha_i + \alpha_{\nu +1 -i} = n-1.$$

V.~V.~Goryunov \cite{Goryunov81} computed the spectra of the simple, uni-
and bimodal hypersurface singularities. Steenbrink
\cite{Steenbrink99} compiled tables of the spectra for all ${\cal
K}$-unimodal ICIS. If $Y=\CC^{n+1}$ and $f$ is Newton non-degenerate then
the spectrum can be computed from the Newton diagram, see \cite{MSaito88,
VKh85}. Other methods to compute the spectrum in the case when $(Y,0)$ is
smooth have been given by Schulze and Steenbrink \cite{SchSt00}.

The most famous property of the spectrum is the semicontinuity conjectured
by Arnold \cite{Arnold81} and proved by Steenbrink \cite{Steenbrink85} for
the case when $(Y,0)$ smooth and the author and Steenbrink in the general case
\cite{ESt98}:

\begin{theorem}[Semicontinuity theorem]
The spectrum behaves semicontinuously under deformation of the singularity
in the following sense: If $f'$ (with $\nu' < \nu$) appears in the
semi-universal deformation of $f$, then
$$\alpha_i \leq \alpha'_i.$$
\end{theorem}

The variance of the spectrum measures the distribution of the spectral
numbers with respect to the central point and is defined by
$$V= \frac{1}{\nu} \sum_{i=1}^\nu \left( \alpha_i - \frac{n-1}{2}
\right)^2.$$

C.~Hertling \cite{Hertling00} proposed the following conjecture

\begin{conjecture}[Hertling]
If $(Y,0)$ is smooth (so $\nu=\mu$), then
$$V \leq \frac{\alpha_\mu - \alpha_1}{12}.$$
\end{conjecture}

One has equality if $f$ is weighted homogeneous, as shown by A.~Dimca
\cite{Dimca00} and Hertling \cite{Hertling00}. M.~Saito \cite{MSaito00}
showed that Hertling's conjecture holds for all irreducible plane curve
singularities. Th.~Br\'elivet \cite{Brelivet02, Brelivet04} showed that
the conjecture holds for all curve singularities. Br\'elivet and Hertling
have stated more general conjectures involving higher moments \cite{BH04}.

Let $Y=\CC^{n+1}$. We shall now give several different interpretations of
the smallest exponent $\alpha_1$.

Let $\omega$ be a holomorphic $(n+1)$-form on $\CC^{n+1}$. For $0 < | t | < \delta$
let $\eta(t)$ be
a continuously varying homology class of dimension $n$ on $X_t$ and
consider the function
$$I(t) = \int_{\eta(t)} \frac{\omega}{df}.$$
This function admits an asymptotic expansion as $t$ tends to zero:
$$I(t)= \sum_{\alpha, q} \frac{1}{q!} C_{\alpha , q}^{\omega , \eta}
t^\alpha (\log t)^q,$$
such that $q \in \ZZ$, $0 \leq q \leq n$, $\alpha \in \QQ$, $\alpha > -1$
and $e^{2 \pi i \alpha}$ is an eigenvalue of the semisimple part of the
monodromy operator. By \cite{Varchenko80} we have
$$\alpha_1= \beta_\CC -1:= \min \{ \alpha \, | \, \exists \omega, \eta , q
\ C_{\alpha , q}^{\omega , \eta} \neq 0\}.$$
The number $\beta_\CC$ is the {\em complex singularity index} (cf.\
\cite{Arnold73}, where in fact the number $\frac{n+1}{2}-\beta_\CC$ is called the complex
singularity index). For a simple singularity in
$\CC^3$, one has $\beta_\CC = 1+ \frac{1}{N}$ where $N$ is the
Coxeter number of the
singularity (cf.\ \cite{Arnold73}).

With the notations of \S 3, let $k_i$ be the multiplicity of $\pi^\ast(dx_0
\wedge \ldots \wedge dx_n)$ along the divisor $E_i$, $i=1, \ldots ,s$. 
Let $m_0$ and $k_0$ be the order of $\widetilde{f}$ and the multiplicity of 
$\pi^\ast(dx_0 \wedge \ldots \wedge dx_n)$ respectively along the divisor $\widetilde{X}_0$. So
$m_0=1$ and $k_0=0$. Let
$$\lambda:= \min \left. \left\{ \frac{k_i+1}{m_i} \, \right| \, i=0, \ldots
,s \right\}.$$
K.-Ch.~Lo \cite{Lo79} has shown that
$$\beta_\CC \geq \lambda.$$
T.~Yano \cite{Yano83} and B.~Lichtin \cite{Lichtin89} have shown that if $\lambda < 1$ then
$$\beta_\CC = \lambda.$$

Ehlers and Lo \cite{EL82} have shown that for a Newton non-degenerate
function,
$\beta_\CC=1/t_0$ where $(t_0, \ldots ,t_0)$ is the intersection point of
the diagonal $t \mapsto (t, \ldots ,t)$ with the Newton diagram of $f$.

J.~Kollar \cite{Kollar97} has shown that $\lambda$ is equal to the {\em log canonical threshold}. 

Moreover, we have the following
relations which were recently brought back into attention in the framework
of multiplier ideals (see e.g.\ \cite{ELSV04}).

For a rational number $\alpha$ we define the following ideal ${\cal
A}_\alpha$ in the ring ${\cal O}_{Y,0}$ of analytic functions on $(Y,0)$:
$${\cal A}_\alpha := \left\{ \phi \in {\cal O}_{Y,0} \, \left| \, \inf_{1 \leq i \leq s}
\left( \frac{1+k_i+ \nu_i(\phi)}{m_i} -1 \right) > \alpha \right\} \right.$$
where $\nu_i(\phi)$ denotes the order of $\phi$ along the divisor $E_i$. This
is a {\em multiplier ideal} in the sense of Y.-T.~Siu and A.~Nadel (see
\cite{ELSV04}).

\begin{definition} We define a sequence of numbers
$$\xi_0=0 < \xi_1 < \xi_2 < \ldots$$
as follows: ${\cal A}_\alpha= {\cal A}_{\xi_i}$ for $\alpha \in [\xi_i,
\xi_{i+1})$ and ${\cal A}_{\xi_{i+1}} \neq {\cal A}_{\xi_i}$ for $i=0,1,
\ldots$. These numbers are called {\em jumping numbers}.
\end{definition}

These numbers first appeared implicitly in a paper of A.~Libgober
\cite{Libgober83}. The above definition is due to Loeser and M.~Vaqui\'e
\cite{LV90, Vaquie92}.

Varchenko \cite{Varchenko82} (see also \cite{Budur03}) proved the following
statement:
$$\alpha \in (-1,0], \quad \alpha \in {\rm Sp}(f) \Leftrightarrow \alpha +1
= \xi_i \mbox{ for some } i.$$
M.~Saito \cite{MSaito93} showed that the Bernstein-Sato polynomial $b_f(s)$
(see \S 1) has roots in
$[0,1)$ which do not come from the spectrum of $f$.


\section{Monodromy and the topology of the singularity}
Let $B_\eps$ be a closed ball as in \S 1 and let $K:=f^{-1}(0)
\cap
\partial B_\eps \cap Y$ be the link of the singularity $(X_0,0)$.

First assume that $Y=\CC^{n+1}$. Milnor \cite{Milnor68} has shown that
the manifold $K$ is a homology sphere (and when $n \neq 2$
actually a topological sphere) if and only if the integer
$$\Delta(1)=\det({\rm id}_\ast - h_\ast)$$
is equal to $\pm 1$.
Let $n \neq 2$ and assume that $K$ is a topological sphere. The
differentiable structure of $K$ is completely determined by
the Kervaire invariant $c(X_\delta) \in \ZZ_2$ if $n$ is odd, or by the
signature of the intersection matrix $A$ if $n$ is even (cf.~\cite{Milnor68}). If $n$ is
odd, then by a theorem of J.~Levine \cite{Levine66} the Kervaire invariant
is given by
$$ c(X_\delta) = \left\{ \begin{array}{cl} 0 & \mbox{if } \Delta(-1)
\equiv \pm 1\, (\mbox{mod }8),\\ 1 & \mbox{if } \Delta(-1) \equiv \pm 3\,
(\mbox{mod }8). \end{array} \right. $$
If $n$ is even and $f$ is weighted homogeneous, then the signature of the 
intersection matrix $A$ is determined by the eigenvalues of the monodromy, see \cite{Steenbrink77a}.
Hence in many cases the complex monodromy operator determines the differentiable structure of $K$.

In fact it is shown in \cite{Milnor68} that $K$ is
$(n-2)$-connected, that the rank of $H_{n-1}(K)$ is equal to the
dimension of the eigenspace of $h_\ast$ corresponding to the eigenvalue 1,
and that, if the rank is equal to zero, the order of $H_{n-1}(K)$  is
equal to $|\Delta(1)|$. Here we use reduced homology if $n=1$.
If $f$ is weighted homogeneous, a formula for the rank of
$H_{n-1}(K)$ and for $\Delta(1)$ in terms of weights and degree is given in
\cite{MO70}.

If $n=1$, Durfee \cite{Durfee75a} relates the topology of a branched cyclic
cover of the link $K$ to
the characteristic polynomial of the monodromy. B.~G.~Cooper
\cite{Cooper82} has calculated the
homology of the link $K$ for some special weighted homogeneous polynomials $f$.

Let $C$ be the matrix of $h_\ast$ with respect to a basis of $H_n(X_\delta,
\CC)$ and let $I$ be the $\mu \times \mu$ identity matrix. In the case
when $f$ is weighted homogeneous, Orlik \cite{Orlik72} has stated the
following conjecture:

\begin{conjecture}[Orlik]
The matrix $C$ can be diagonalized over the integers, i.e.\ there exist
unimodular matrices $U(t)$ and $V(t)$ with entries in the ring $\ZZ[t]$ so
that
$$U(t)(tI-C)V(t) = {\rm diag}(m_1(t), \ldots , m_\mu(t))$$
where $m_i(t)$ divides $m_{i+1}(t)$ for $i=1, \ldots ,\mu -1$.
\end{conjecture}

Since the ring $\CC[t]$ is a principal ideal domain, such matrices exist
over $\CC[t]$. The conjecture implies that
$$H_{n-1}(K) = \ZZ_{m_1(1)} \oplus \ldots \oplus \ZZ_{m_\mu(1)}$$
where $\ZZ_1$ is the trivial group and $\ZZ_0$ is the infinite cyclic
group.

The conjecture holds for $f$ weighted homogeneous and $n=2$ as follows from
\cite{OW71}.

Sometimes the conjecture is extended to germs $f$ with finite monodromy. Then the
following is known about the more general conjecture. From \cite{A'Campo73b} one can
derive that Orlik's conjecture is true for irreducible plane curve
singularities. F.~Michel and C.~Weber \cite{MW84, MW86} have shown that
Orlik's conjecture is false for plane curve singularities with more than
one branch.

Let $n=2$ and $f$ be weighted homogeneous with weights $q_1$, $q_2$, $q_3$
and degree $d$. Then Y.~Xu and S.-T.~Yau \cite{XY89} have shown that the
characteristic polynomial $\Delta(t)$ of the monodromy and the fundamental
group $\pi_1(K)$ of the link determine the embedded topological type of
$(X_0,0)$. Let $K$ be in addition a rational homology sphere. Then
R.~Mendris and A.~N\'emethi \cite{MN02} have observed that it follows from
\cite{Ebeling02} that $\Delta(t)$ is already determined by $\pi_1(K)$.
Define
$$R:=d-q_1-q_2-q_3.$$
N\'emethi and L.~I.~Nicolaescu \cite{NN04}
have derived from \cite{Ebeling02} that
$$\frac{\Delta(t)}{\Delta(1)} = 1 +
\frac{\mu}{2}(t-1) + \ldots \quad \mbox{and} \quad
\frac{\Delta^\ast(t)}{\Delta^\ast(1)} = 1 + \frac{R}{2}(1-e_{\rm st})(t-1) +
\ldots$$
where $e_{\rm st}$ is  Batyrev's stringy Euler characteristic of
$(X_0,0)$ (cf.\ \cite{Batyrev98}) as generalized by Veys in \cite{Veys04}.

Now let $Y=\CC^{n+1}$ and let $f$ be again general. The matrix $V$ of \S 2
is the matrix of the (integral) Seifert form of the singularity
$(X_0,0)$. If $n \geq 3$ then results of M.~Kervaire \cite{Kervaire65} and
J.~Levine \cite{Levine70} show that the Seifert form determines the
(embedded) topological type of the singularity, see also \cite{Durfee74}.
If $n=1$ and $f$ defines
an irreducible curve singularity, then it follows from \cite{Burau32} and
\cite{Zariski32} that the
integral monodromy and even the rational monodromy determines the topology
of the singularity.
M.-C.~Grima \cite{Grima74} has given examples of plane curve singularities with
two branches of different topological types with the same rational
monodromy, but different integral
monodromy.
Ph.~Du
Bois and Michel
\cite{DBM91, DBM94} have shown that the integral Seifert form does not
always determine the topology
of the singularity in the case
$n=1$. Using suspensions of the examples of Du Bois and
Michel, Artal-Bartolo \cite{ArtalBartolo91b} has shown that the same applies
to the case $n=2$.

Let $(Y,0)$ be a weighted homogeneous ICIS and let $f$ be weighted
homogeneous. Let $L:=\partial B_\eps \cap Y$. Dimca \cite{Dimca85} has shown
that for $n \geq 2$ one has the following formula for the Betti numbers of
$L$ and $K$:
$$b_{n+1}(L) + b_n(K) = {\rm dim} \, \ker ({\rm Id}_\ast - h_\ast).$$
Hamm \cite{Hamm72} computed the characteristic polynomial of the monodromy
for some ICIS which are
generalizations of the Brieskorn-Pham singularities, the
Brieskorn-Hamm-Pham singularities.  The
homology torsion of the link of a Brieskorn-Hamm-Pham singularity was
computed by Randell
\cite{Randell75}.

%
%
\section*{Acknowledgements}
This research was partially supported by the DFG-programme ''Global methods in
complex geometry'' (Eb 102/4--3).
The author is grateful to Th.~Eckl for drawing his attention to jumping
numbers of multiplier ideals
and the connection with the spectrum (see \S 4). He would like to thank A.~Melle-Hern\'andez and
W.~Veys for their useful comments.
%
%

%


\begin{thebibliography}{999}
%
%
%
%



\bibitem{A'Campo73a}  A'Campo, N.: {\it Le nombre de Lefschetz d'une
monodromie.} Indag. Math. {\bf 35}, 113--118 (1973).

\bibitem{A'Campo73b}  A'Campo, N.: {\it Sur la monodromie des
singularit\'{e}s isol\'{e}es d'hypersurfaces complexes.} Invent. math. {\bf
20}, 147--169 (1973).

\bibitem{A'Campo75a}  A'Campo, N.: {\it La fonction z\^eta d'une
monodromie.}
Comment. Math. Helv. {\bf 50}, 233--248 (1975).

\bibitem{A'Campo75b}  A'Campo, N.: {\it Le groupe de monodromie du
d\'eploiement des singularit\'es isol\'ees de courbes planes I.} Math. Ann.
{\bf 213}, 1--32 (1975).

\bibitem{A'Campo75c}  A'Campo, N.:
{\it Le groupe de monodromie du d\'eploiement des singularit\'es isol\'ees de
courbes planes. II.}
In: Act. Cong. Int. Math. (Vancouver 1974),
Canad. Math. Congress, Montreal 1975, Vol. 1, pp. 395--404.

\bibitem{A'Campo03}  A'Campo, N.: {\it Monodromy of real isolated
singularities.}
Topology {\bf 42}, 1229--1240 (2003).

\bibitem{Arnold73}  Arnold, V.~I.: {\it Remarks on the method of
stationary phase and on the
Coxeter numbers.} Uspehi Mat. Nauk {\bf 28}:5, 17--44 (1973) (Engl.
translation in Russian Math.
Surveys {\bf 28}:5, 19--48 (1973)).

\bibitem{Arnold75}  Arnold, V.~I.: {\it Critical points of smooth
functions.} In: Proc. Intern Congr. Math. Vancouver, 1974, Canad. Math.
Congress, Montreal 1975, Vol.
1, pp. 19--39.

\bibitem{Arnold76}  Arnold, V.~I.: {\it Local normal forms of functions.}
Invent.
math. {\bf 35}, 87--109 (1976).

\bibitem{Arnold81}  Arnold, V.~I.: {\it On some problems in singularity
theory.} In: {\it Geometry and Analysis}, Papers dedicated to the memory of
V.~K.~Patodi, Bombay, 1981, pp. 1--10.

\bibitem{AGV88} Arnold,  V.~I.,  Gusein-Zade, S.~M.,  Varchenko, A.~N.:
{\it Singularities of Differentiable Maps}, Volume II. Birkh\"auser, Boston Basel
Berlin 1988.

\bibitem{ArtalBartolo91a}  Artal-Bartolo, E.: {\it Sur la
monodromie des singularit\'es superisol\'ees.} C. R. Acad. Sci. Paris S\'er.
I Math. {\bf 312}, 601--604 (1991).

\bibitem{ArtalBartolo91b}  Artal-Bartolo, E.: {\it Forme de
Seifert des singularit\'es
de surface.} C. R. Acad. Sci. Paris S\'er. I Math. {\bf 313}, 689--692 (1991).

\bibitem{ArtalBartolo94} Artal-Bartolo, E.: {\it Forme de Jordan
de la monodromie des singularit\'es superisol\'ees de surfaces.} Mem. Amer.
Math. Soc. {\bf 109}, no. 525, (1994).

\bibitem{ACLM02} Artal-Bartolo, E., Cassou-Nogu\`es, P., Luengo, I.,
Melle-Hern\'andez, A.: {\it Monodromy conjecture for some surface singularities.}
Ann. Scient. Ec. Norm. Sup. {\bf 35}, 605--640 (2002).

\bibitem{ACLM03} Artal-Bartolo, E. Cassou-Nogu\`es, P., Luengo, I.,
Melle-Hern\'andez, A.: {\it Quasi-ordinary power series and their zeta functions.}
Preprint 2003, Memoirs Amer. Math. Soc. (to appear).

\bibitem{ALM05} Artal-Bartolo, E., Luengo, I., Melle-Hern\'andez, A.: {\it Superisolated
singularities.} These proceedings.

\bibitem{Barlet86}  Barlet, D.:  {\it Monodromie et p\^oles du
prolongement m\' eromorphe de
$\int\sb X\vert f\vert
\sp {2\lambda}\square$.} Bull. Soc. Math. France {\bf 114}, 247--269 (1986).

\bibitem{Batyrev98} Batyrev, V.~V.: {\it Stringy Hodge numbers of
varieties with Gorenstein
canonical singularities.} In: {\it Integrable systems and algebraic geometry}
(Kobe/Kyoto, 1997),
World Sci. Publishing, River Edge, NJ, 1998, pp. 1--32.

\bibitem{Bey72a} Bey, C.~H.:  {\it Sur l'irr\'eductibilit\'e de la
monodromie locale.}
C. R. Acad. Sci. Paris S\'er. A-B {\bf 275}, A21--A24 (1972).

\bibitem{Bey72b} Bey, C.~H.: {\it Sur l'irr\'eductibilit\'e de la
monodromie locale; application
\`a l'\'equisingularit\'e.}  C. R. Acad. Sci. Paris S\'er. A-B {\bf 275},
A105--A107 (1972).

\bibitem{Bony76} Bony, J.-M.: {\it Polyn\^omes de Bernstein et monodromie
(d'apr\`es B.
Malgrange).}  S\'eminaire Bourbaki (1974/1975), Exp. No. 459, Lecture Notes in
Math., Vol. {\bf 514}, Springer, Berlin, 1976, pp. 97--110.

\bibitem{Bourbaki68} Bourbaki, N.: {\it Groupes et alg\`{e}bres de
Lie,} Chapitres 4,5 et 6.
Hermann, Paris, 1968.

\bibitem{Brelivet02} Br\'elivet, Th.: {\it Variance of spectral numbers
and Newton polygons.} Bull. Sci. Math. {\bf 126}, 333--342 (2002).

\bibitem{Brelivet04} Br\'elivet,  Th.: {\it The Hertling
conjecture in dimension 2.} Preprint, math.AG/0405489.

\bibitem{BH04} Br\'elivet, Th., Hertling, C.: {\it Bernoulli moments of
spectral numbers and Hodge numbers.} Preprint, math.AG/0405501.

\bibitem{BS74} Brian\c{c}on, J., Skoda, H.: {\it Sur la cl\^{o}ture
int\'egrale d'un id\'eale de germes de fonctions holomorphes en un point de
$\CC^n$.} C. R. Acad. Sci. Paris S\'er. I Math. {\bf 278}, 949 (1974).

\bibitem{Brieskorn70} Brieskorn, E.: {\it Die Monodromie der
isolierten Singularit\"aten von Hyperfl\"achen.} Manuscripta math. {\bf 2},
103--160 (1970).

\bibitem{Brieskorn73} Brieskorn, E.: {\it Vue d'ensemble sur les
probl\`emes de monodromie.} In: {\it Singularit\'es \`a Carg\`ese}, Ast\'erisque
{\bf 7} et {\bf 8}, 195--212 (1973).

\bibitem{Brieskorn76} Brieskorn, E.: {\it Singularit\"aten.} Jber.
DMV {\bf 78}, 93--112 (1976).

\bibitem{Budur03} Budur, N.: {\it On Hodge spectrum and multiplier
ideals.} Math. Ann. {\bf 327},
257--270 (2003).

\bibitem{Burau32} Burau, W.: {\it Kennzeichnung der Schlauchknoten.} Abh.
Math. Sem. Hamburg {\bf
9}, 125--133 (1932).

\bibitem{Clemens69} Clemens, C.~H.: {\it Picard-Lefschetz theorem for
families of nonsingular
algebraic varieties acquiring ordinary singularities.} Trans. Amer. Math.
Soc. {\bf 136}, 93--108
(1969).

\bibitem{Cooper82} Cooper, B.~G.: {\it On the monodromy at isolated
singularities of weighted
homogeneous polynomials.}  Trans. Amer. Math. Soc. {\bf 269}, 149--166 (1982).

\bibitem{Deligne70} Deligne, P.: {\it \'Equations Diff\'erentielles \`a
Points Singuliers R\'eguliers.} Lecture Notes in Math., Vol. {\bf 163},
Springer-Verlag, Berlin etc., 1970.

\bibitem{DK73} Deligne, P., Katz, N.: {\it Groupes de monodromie en
g\'eom\'etrie
alg\'ebrique II} (SGA 7 II). Lecture Notes in Math., Vol. {\bf 340},
Springer-Verlag, Berlin New York,
1973.

\bibitem{Denef91} Denef, J.: {\it Report on Igusa's local zeta function.},
S\'eminaire Bourbaki, Ast\'erisque
{\bf 201-202-203}, 358--386 (1991).

\bibitem{DL92} Denef, J., Loeser, F.: {\it Charact\'eristiques
d'Euler-Poincar\'e, fonctions z\^eta locales, et modifications analytiques.}
J. Amer. Math. Soc. {\bf 5}, 705--720 (1992).

\bibitem{Dimca85} Dimca, A.: {\it Monodromy and Betti numbers of
weighted complete intersections.}
Topology {\bf 24}, 369--374 (1985).

\bibitem{Dimca00} Dimca, A.: {\it Monodromy and Hodge theory of regular
functions.} In: {\it New Developments in
Singularity Theory} (Siersma, D., et al., eds.), NATO Science Series II, Vol.
{\bf 21}, Kluwer Acad. Publ., Dordrecht et al., 2001, pp. 257--278.

\bibitem{vDSt89} van Doorn, M.~G.~M., Steenbrink, J.~H.~M.: {\it  A
supplement to the monodromy theorem.} Abh. Math. Sem. Univ. Hamburg {\bf
59}, 225--233 (1989).

\bibitem{DBM91} Du Bois, Ph., Michel, F.: {\it  Sur la forme de Seifert des
entrelacs alg\'ebriques.}
C. R. Acad. Sci. Paris S\'er. I Math. {\bf 313}, 297--300 (1991).

\bibitem{DBM92} Du Bois, Ph., Michel, F.: {\it  Filtration par le poids et
monodromie enti\`ere.}
Bull. Soc. Math. France {\bf 120}, 129--167 (1992).

\bibitem{DBM94} Du Bois, Ph., Michel, F.: {\it  The integral Seifert form
does not determine the
topology of plane curve germs.} J. Algebraic Geom. {\bf 3}, 1--38 (1994).

\bibitem{Durfee74} Durfee, A.~H.: {\it  Fibered knots and algebraic
singularities.} Topology {\bf
13}, 47--59 (1974).

\bibitem{Durfee75a} Durfee, A.~H.: {\it  The characteristic polynomial
of the monodromy.}
Pacific J. Math. {\bf 59}, 21--26 (1975).

\bibitem{Durfee75b} Durfee, A.~H.: {\it  The monodromy of a
degenerating family of curves.}
Invent. Math. {\bf 28}, 231--241 (1975).

\bibitem{Ebeling80} Ebeling, W.: {\it  Quadratische Formen und
Monodromiegruppen von Singularit\"aten.} Math.\ Ann.\ {\bf 255}, 463--498 (1981).

\bibitem{Ebeling83} Ebeling, W.: {\it  Milnor lattices and geometric
bases of some special singularities.}
          In: {\it  N{\oe}uds, tresses et singularit\'es} (Ed. Weber, C.),
          Monographie Enseign. Math. {\bf 31},
          Gen{\`e}ve, 1983, pp. 129--146;
          Enseign. Math. {\bf 29}, 263--280 (1983).

\bibitem{Ebeling85} Ebeling, W.: {\it  The Milnor lattices of the
elliptic hypersurface singularities.}
          Proc.\ London Math.\ Soc.\ (3), {\bf	53}, 85--111 (1986).

\bibitem{Ebeling87} Ebeling, W.: {\it  The Monodromy Groups of Isolated
Singularities of Complete Intersections.}	 Lecture Notes in Math., Vol. {\bf
1293}, Springer-Verlag, Berlin etc., 1987.

\bibitem{Ebeling96} Ebeling, W.: {\it  On Coxeter-Dynkin diagrams of
hypersurface singularities.}
J. Math. Sci. {\bf 82}, 3657--3664 (1996).

\bibitem{Ebeling99} Ebeling, W.: {\it  Strange duality, mirror
symmetry, and the
Leech lattice.} In: {\it  Singularity Theory}, Proceedings of the European
Singularities Conference, Liverpool 1996 (Bruce, J.~W., Mond, D., eds.), London
Math. Soc. Lecture Note Ser., Vol. {\bf 263}, Cambridge University Press,
Cambridge, 1999, pp. 55--77.

\bibitem{Ebeling00} Ebeling, W.: {\it  Strange duality and polar duality.}
 J. London Math. Soc. {\bf 61}, 823--834 (2000).

\bibitem{Ebeling01} Ebeling, W.: {\it  Funktionentheorie,
Differentialtopologie und
Singularit\"aten.} Vieweg, Braunschweig, Wiesbaden, 2001.

\bibitem{Ebeling02} Ebeling, W.: {\it  Poincar\'{e} series and
monodromy of a
two-dimensional quasihomogeneous hypersurface singularity.}
Manuscripta math. {\bf 107}, 271--282 (2002).

\bibitem{Ebeling05} Ebeling, W.: {\it  Mirror symmetry, Kobayashi's
duality, and Saito's duality.}
Preprint 2005.

\bibitem{EG95} Ebeling, W., Gusein-Zade, S.~M.: {\it   Coxeter-Dynkin
diagrams of the complete
intersection singularities $Z_9$ and $Z_{10}$.} Math. Z. {\bf 218}, 549--562
(1995).

\bibitem{EG96} Ebeling, W., Gusein-Zade, S.~M.: {\it   Coxeter-Dynkin
diagrams of fat points in
$\CC^2$ and of their stabilizations.} Math.\ Ann.\ {\bf 306}, 487--512 (1996).

\bibitem{EG98} Ebeling, W., Gusein-Zade, S.~M.: {\it  Suspensions of fat points
and their intersection forms.} In: {\it Singularities.} The Brieskorn Anniversary
Volume (Arnold, V.~I., Greuel, G.-M., Steenbrink, J.~H.~M., eds.), Progr. Math.,
Vol. {\bf 162}, Birkh\"auser, Basel 1998, pp. 141--165.

\bibitem{EG99} Ebeling, W., Gusein-Zade, S.~M.: {\it  On the index of a vector
field at an isolated singularity.} In: {\it  The Arnoldfest}: Proceedings of a
Conference in Honour of V.~I.~Arnold for his Sixtieth Birthday (Bierstone, E.,
Khesin, B., Khovanskii, A., Marsden, J., eds.), Fields Institute Communications,
Vol. {\bf 24}, Amer. Math. Soc., Providence 1999, pp. 141--152.

\bibitem{EG02} Ebeling, W., Gusein-Zade, S.~M.: {\it  Poincar\'e series and
zeta function of the monodromy of a quasihomogeneous singularity.} Math.
Res. Lett. {\bf 9}, 509--513 (2002).

\bibitem{EG04} Ebeling, W.,
Gusein-Zade, S.~M.: {\it  Monodromies and Poincar\'e series of quasihomogeneous
complete intersections.} Abh. Math. Sem. Univ. Hamburg {\bf 74}, 175--179
(2004).

\bibitem{ESt98} Ebeling, W., J.~H.~M. Steenbrink: {\it   Spectral pairs for
isolated complete
intersection singularities.} J. Alg. Geom. {\bf 7}, 55--76 (1998).

\bibitem{EW85} Ebeling, W., C.~T.~C. Wall: {\it  Kodaira singularities and
an extension of
Arnold's strange duality.} Compositio Math.\ {\bf 56}, 3--77 (1985).

\bibitem{Ehlers78} Ehlers, F.: {\it  Newtonpolyeder und die Monodromie
von Hyper\-fl\"achen\-singularit\"aten.} Bonner Mathematische Schriften {\bf
111}, Bonn 1978.

\bibitem{EL82} Ehlers, F., Lo, K.-Ch.: {\it  Minimal characteristic exponent
of the Gauss-Manin connection of isolated singular point and Newton
Polyhedron.} Math. Ann. {\bf 259}, 431--441 (1982).

\bibitem{EN85} Eisenbud, D., Neumann, W.: {\it  Three-dimensional link theory
and invariants of plane
curve singularities.} Annals of Math. Studies Vol. {\bf 110}, Princeton
University Press,
Princeton, 1985.

\bibitem{ELSV04} Ein, L., Lazarsfeld, R., Smith, K., Varolin, D.: {\it 
Jumping coefficients of multiplier ideals.} Duke Math. J. {\bf 123},
469--506 (2004).

\bibitem{Gabrielov73} Gabrielov, A.~M.: {\it  Intersection matrices
for certain
singularities.} Funktsional. Anal. i Prilozhen. {\bf 7}:3, 18--32 (1973)
(Engl. translation in Funct. Anal. Appl. {\bf 7}, 182--193 (1974)).

\bibitem{Gabrielov74a} Gabrielov, A.~M.: {\it   Bifurcations, Dynkin
diagrams and the modality
of isolated singularities.} Funktsional. Anal. i Prilozhen. {\bf 8}:2,
7--12 (1974). (Engl.
translation: Funct. Anal. Appl. {\bf 8}, 94--98 (1974).)

\bibitem{Gabrielov74b} Gabrielov, A.~M.: {\it  Dynkin diagrams of
unimodal
singularities.} Funktsional. Anal. i Prilozhen. {\bf 8}:3, 1--6 (1974) (Engl.
translation in Funct. Anal. Appl. {\bf 8}, 192--196 (1974)).

\bibitem{Gabrielov79} Gabrielov, A.~M.: {\it 
Polar curves and intersection matrices of singularities.} Invent. math. {\bf
54}, 15--22 (1979).

\bibitem{Giusti77} Giusti, M.: {\it  Classification des singularit\'es
isol\'ees
d'intersections compl\`etes simples.} C.R. Acad. Sc. Paris, S\'er. A, {\bf
284}, 167--170 (1977).

\bibitem{Goryunov81} Goryunov, V.~V.: {\it  Adjacencies of spectra of
certain singularities.} Vestnik MGU Ser. Math. {\bf 4}, 19--22 (1981) (in
Russian).

\bibitem{Greuel75} Greuel, G.-M.: {\it  Der Gau\ss-Manin-Zusammenhang
isolierter Singularit\"aten von vollst\"andigen Durchschnitten.}
Math. Ann. {\bf 214}, 235--266 (1975).

\bibitem{GH78} Greuel, G.-M., Hamm, H.: {\it  Invarianten quasihomogener
vollst\"andiger
Durchschnitte.} Invent. math. {\bf 49}, 67--86 (1978).

\bibitem{Grima74} Grima, M.-C.: {\it  La monodromie rationnelle ne
d\'etermine pas la topologie d'une
hypersurface complexe.} In: {\it  Fonctions de plusieurs variables complexes} (S\'em.
Fran\c{c}ois Norguet,
1970--1973; \`a la m\'emoire d'Andr\'e Martineau), Lecture Notes in Math., Vol.
{\bf 409},
Springer, Berlin, 1974, pp. 580--602.

\bibitem{Griffiths70} Griffiths, P.~A.: {\it  Periods of integrals on
algebraic manifolds:
Summary of main results and discussion of open problems.} Bull. Amer. Math.
Soc. {\bf 76}, 228--296
(1970).

\bibitem{Grothendieck72} Grothendieck, A.: {\it  Groupes de
monodromie en g\'eom\'etrie
alg\'ebrique I} (SGA 7 I). Lecture Notes in Math., Vol. {\bf 288},
Springer-Verlag, Berlin New York,
1972.

\bibitem{GuseinZade74a} Gusein-Zade, S.~M.: {\it  Intersection
matrices for certain singularities of functions of two variables.}
Funktsional. Anal. i Prilozhen. {\bf 8}, 11--15 (1974) (Engl. translation in
Funct. Anal. Appl. {\bf 8}, 10--13 (1974)).

\bibitem{GuseinZade74b} Gusein-Zade, S.~M.: {\it  Dynkin
diagrams for singularities of functions of two variables.} Funktsional. Anal. i Prilozhen.
 {\bf 8}, 23--30 (1974) (Engl. translation in Funct.
Anal. Appl. {\bf 8}, 295--300 (1974)).

\bibitem{GuseinZade76} Gusein-Zade, S.~M.: {\it  Characteristic
polynomial of classical
monodromy for singular series.} Funktsional. Anal. i Prilozhen. {\bf 10},
78--79 (1976)
(Engl. translation in Funct. Anal. Appl. {\bf 10}, 229--231 (1976)).

\bibitem{GuseinZade84} Gusein-Zade, S.~M.: {\it  The index of a
singular point of a gradient
vector field.} Funktsional. Anal. i Prilozhen. {\bf 18}, 7--12 (1984) (Engl.
translation in
Funct. Anal. Appl. {\bf 18}, 6--10 (1984)).

\bibitem{CDG99} Gusein-Zade, S.~M., Delgado, F., Campillo, A.: {\it On the
monodromy of
a plane curve singularity and the Poincar\'{e} series of the ring of
functions on the
curve.} Funktsional. Anal. i Prilozhen. {\bf 33}:1, 66--68 (1999) (Engl.
translation in Funct. Anal. Appl. {\bf 33}:1, 56--57 (1999)).

\bibitem{GLM97} Gusein-Zade, S.~M., Luengo, I., Melle-Hern\'andez, A.: 
{\it  Partial resolutions and the zeta-function
of a singularity.}
Comment. Math. Helv. {\bf 72}, 244--256 (1997).

\bibitem{GLM01} Gusein-Zade, S.~M., Luengo, I., Melle-Hern\'andez, A.: {\it 
Bifurcations and topology of meromorphic germs.}
In: {\it  New Developments in
Singularity Theory} (Siersma, D., et al., eds.), NATO Science Series II, Vol.
{\bf 21}, Kluwer Acad. Publ., Dordrecht et al., 2001, pp. 279--304.

\bibitem{GLM04} Gusein-Zade, S.~M., Luengo, I., Melle-Hern\'andez, A.: {\it 
Integration over spaces of non-parametrized arcs and motivic versions of
the monodromy zeta function.} Preprint 2004.

\bibitem{Hamm71} Hamm,  H.: {\it  Lokale topologische Eigenschaften komplexer
R\"aume.} Math. Ann. {\bf 191}, 235--252 (1971).

\bibitem{Hamm72} Hamm, H.: {\it  Exotische Sph\"aren als Umgebungsr\"ander
in speziellen komplexen
R\"aumen.}  Math. Ann. {\bf 197}, 44--56 (1972).

\bibitem{HL74} Hefez, A., Lazzeri, F.: {\it  The intersection matrix of
Brieskorn singularities.}
Invent. Math. {\bf 25}, 143--157 (1974).

\bibitem{Hertling00} Hertling, C.: {\it  Frobenius manifolds and
variance of the spectral numbers.} In: {\it  New Developments in
Singularity Theory} (Siersma, D., et al., eds.), NATO Science Series II, Vol.
{\bf 21}, Kluwer Acad. Publ., Dordrecht et al., 2001, pp. 235--255.

\bibitem{Igusa75} Igusa, J.: {\it  Complex powers and asymptotic
expansions I,} J. Reine Angew.
Math. {\bf 268/269}, 110--130 (1974); II, ibid. {\bf 278/279}, 307--321 (1975).

\bibitem{Il'yuta87} Il'yuta, G.~G.: {\it 
On the Coxeter transformation of an isolated singularity.}
Uspekhi Mat. Nauk {\bf 42}:2, 227--228 (1987)
(Engl. translation in Russian Math. Surveys {\bf 42}:2, 279--280 (1987)).

\bibitem{Katz70} Katz, N.~M.: {\it  Nilpotent connections and the
monodromy theorem: Applications
of a result of Turrittin.}  Inst. Hautes \'Etudes Sci. Publ. Math.  {\bf 39},
175--232 (1970).

\bibitem{Kervaire65} Kervaire, M.: {\it  Les n{\oe}uds de dimension
sup\'erieures.} Bull. Soc.
Math. France {\bf 93}, 225--271 (1965).

\bibitem{Kollar97} Koll\'ar, J.: {\it  Singularities of pairs.} Proc.
Symp. Pure Math. {\bf 62}, Part 1,
221--287 (1997).

\bibitem{Lamotke75} Lamotke, K.: {\it  Die Homologie isolierter
Singularit\"aten.}
Math. Z. {\bf 143}, 27--44 (1975).

\bibitem{Lamotke05} Lamotke, K.: {\it  Riemannsche Fl\"achen.}
Springer-Verlag, Berlin, Heidelberg, 2005.

\bibitem{Landman67} Landman, A.: {\it  On the Picard-Lefschetz formula
for algebraic manifolds
acquiring general singularities.} Thesis, Berkeley, 1967.

\bibitem{Landman73} Landman, A.: {\it  On the Picard-Lefschetz
transformation for algebraic
manifolds acquiring general singularities.} Trans. Amer. Math. Soc. {\bf
181}, 89--126 (1973).

\bibitem{Lazzeri73} Lazzeri, F.: {\it  A theorem on the monodromy of
isolated singularities.} In: {\it 
Singularit\'es \`a Carg\`ese,} Ast\'erisque
{\bf 7} et {\bf 8}, 269--275 (1973).

\bibitem{Lazzeri88} Lazzeri, F.: {\it  Some remarks on the
Picard-Lefschetz monodromy.} In: {\it 
Introduction \`a la th\'eorie des singularit\'es,} I,
Travaux en Cours, {\bf 36},
Hermann, Paris, 1988, pp. 125--134.

\bibitem{Le72} L\^e D\~ung Tr\'ang: {\it  Sur les n{\oe}uds alg\'ebriques.}
Compos. Math. {\bf 25}, 281--321 (1972).

\bibitem{Le75} L\^e D\~ung Tr\'ang: {\it  La monodromie n'a pas de points fixes.}
J. Fac. Sci. Univ. Tokyo Sect. IA Math. {\bf 22}, 409--427 (1975).

\bibitem{Le78} L\^e D\~ung Tr\'ang: {\it  The geometry of the monodromy theorem.}
In: {\it  C.~P.~Ramanujam, a Tribute} (Ramanathan, K.~G., ed.), Tata Institute
Studies in Math., Vol. {\bf 8}, Springer-Verlag, Berlin etc., 1978.

\bibitem{Lefschetz24} Lefschetz, S.: {\it  L'Analysis situs et la
g\'eom\'etrie alg\'ebrique.}
Gauthiers-Villars, Paris 1924.

\bibitem{Levine66} Levine, J.: {\it  Polynomial invariants of knots of
codimension two.}  Ann. of
Math.  {\bf 84}, 537--554 (1966).

\bibitem{Levine70} Levine, J.: {\it  An algebraic classification of some
knots of codimension
two.} Comment. Math. Helv. {\bf 45}, 185--198 (1970).

\bibitem{Libgober83} Libgober, A.: {\it  Alexander invariants of plane
algebraic curves.} Proc. Symp. Pure Math. Vol.  {\bf 40},
Part 2, 135--143 (1983).

\bibitem{Lichtin89} Lichtin, B.: {\it  Poles of $\vert f(z,w)\vert \sp
{2s}$ and roots of the
$b$-function.} Ark. Mat. {\bf 27}, 283--304 (1989).

\bibitem{Lo79} Lo, K.-Ch.: {\it  Exposants de Gauss-Manin.} In: Pham, F.,
{\it Singularit\'es des Syst\`emes Diff\'erentiels de Gauss-Manin.} Progress in
Mathematics, Vol. {\bf 2}, Birkh\"auser, Basel, Boston, Stuttgart, 1979,
pp. 171--212.

\bibitem{Loeser88} Loeser, F.: {\it  Fonctions d'Igusa $p$-adiques et
polyn\^omes de Bernstein.} Amer. J. Math. {\bf
110}, 1--22 (1988).

\bibitem{Loeser90} Loeser, F.: {\it  Fonctions d'Igusa $p$-adiques,
polyn\^omes de Bernstein, et poly\`edres de Newton.} J. reine angew. Math.
{\bf 412}, 75--96 (1990).

\bibitem{LV90} Loeser, F., Vaqui\'e, M.: {\it  Le polyn\^ome d'Alexander
d'une courbe plane projective.} Topology {\bf 29}, 163--173 (1990).

\bibitem{Lonne03} L\"onne, M.: {\it  Braid monodromy of hypersurface
singularities.} Habilitationsschrift, Hannover, 2003.

\bibitem{Looijenga84} Looijenga, E.~J.~N.: {\it  Isolated Singular
Points on
Complete Intersections.} Cambridge University Press, Cambridge 1984.

\bibitem{Malgrange73}  Malgrange, B.: {\it  Letter to the editors.}
Invent. math. {\bf 20}, 171--172 (1973).

\bibitem{Malgrange74a}  Malgrange, B.: {\it  Sur les polyn\^omes de
I. N. Bernstein.}
S\'eminaire Goulaouic-Schwartz 1973--1974: \'Equations aux d\' eriv\'ees
partielles et analyse
fonctionnelle, Exp. No. 20, 10 pp., Centre de Math., \'Ecole Polytech.,
Paris, 1974.

\bibitem{Malgrange74b}  Malgrange, B.: {\it  Int\'egrales asymptotiques
et monodromie.} Ann. Sc. Ec. Norm. Sup. {\bf 7}, 405--430 (1974).

\bibitem{Malgrange75}  Malgrange, B.: {\it  Le polyn\^ome de Bernstein
d'une singularit\'e
isol\'ee.}  In: {\it  Fourier integral operators and partial differential
equations} (Colloq. Internat.,
Univ. Nice, Nice, 1974),  Lecture Notes in Math., Vol. {\bf 459}, Springer,
Berlin, 1975,
pp. 98--119.

\bibitem{MN02} Mendris, R., N\'emethi, A.: {\it  The link of $\{
f(x,y)+z^n=0\}$ and Zariski's
conjecture.} Compositio Math. {\bf 141}, 502--524 (2005).

\bibitem{MW84} Michel, F., Weber, C.: {\it  Une singularit\'e isol\'ee dont
la monodromie n'admet pas
de forme de Jordan sur les entiers.} C. R. Acad. Sc. Paris S\'er. I Math.
{\bf 299}, 383--386 (1984).

\bibitem{MW86} Michel, F., Weber, C.: {\it  Sur le r\^ole de la monodromie
enti\`ere dans la topologie
des singularit\'es.} Ann. Inst. Fourier {\bf 36}, 183--218 (1986).

\bibitem{Milnor68} Milnor, J.: {\it  Singular Points of Complex
Hypersurfaces.} Ann. of
Math. Studies Vol. {\bf 61}, Princeton University Press, Princeton, 1968.

\bibitem{MO70} Milnor, J., Orlik, P.: {\it  Isolated singularities defined by
weighted
homogeneous polynomials.} Topology {\bf 9}, 385--393 (1970).

\bibitem{NN04} N\'emethi, A., Nicolaescu, L.~I.: {\it  Seiberg-Witten
invariants and surface
singularities. II: Singularities with good $\CC^\ast$-action.} J. London
Math. Soc. {\bf 69},
593--607 (2004).

\bibitem{Oka90} Oka, M.: {\it  Principal zeta-function of non-degenerate
complete
intersection singularity.} J. Fac. Sci. Univ. Tokyo Sect. IA, Math.  {\bf 37},
11--32 (1990).

\bibitem{Orlik72} Orlik, P.: {\it  On the homology of weighted
homogeneous manifolds.} In:
Proceedings of the Second Conference on Compact Transformation Groups
(Univ. Massachusetts, Amherst,
Mass., 1971), Part I, Lecture Notes in Math., Vol. {\bf 298}, Springer,
Berlin, 1972, pp.
260--269.

\bibitem{OR77} Orlik, P., Randell, R.: {\it  The monodromy of weighted
homogeneous singularities.}
Invent. Math. {\bf 39}, 199--211 (1977).

\bibitem{OW71} Orlik, P., Wagreich, Ph.: {\it  Isolated singularities of
algebraic surfaces with
C$\sp{*}$ action.} Ann. of Math. (2) {\bf 93}, 205--228 (1971).

\bibitem{Pham65} Pham, F.: {\it  Formules de Picard-Lefschetz
g\'en\'eralis\'ees et ramification des
int\'egrales.} Bull. Soc. Math. France {\bf 93}, 333--367 (1965).

\bibitem{PS97} Picard, E., Simart, S.: {\it Trait\'e des fonctions
alg\'ebriques de deux variables.}
Vol. I. Gauthier-Villars, Paris, 1897.

\bibitem{Randell75} Randell, R.~C.: {\it  The homology of generalized
Brieskorn manifolds.}
Topology {\bf 14}, 347--355 (1975).

\bibitem{Riemann57} Riemann, B.: {\it  Beitr\"age zur Theorie der durch
die Gau{\ss}sche Reihe $F(\alpha,\beta,\gamma,x)$ darstellbaren Functionen}
(G\"ottingen 1857). In: {\it  Gesammelte mathematische Werke.} Leipzig: Teubner
1892; Springer, Berlin und Teubner, Leipzig 1990, pp. 67-83.

\bibitem{RV03} Rodrigues, B., Veys, W.: {\it  Poles of zeta functions on
normal surfaces.} Proc. London Math. Soc. {\bf 87}, 164--196 (2003).

\bibitem{KSaito88} Saito, K.: {\it  On the existence of exponents prime
to the Coxeter number.} J.
Algebra {\bf 114}, 333--356 (1988).

\bibitem{Saito98a} Saito, K.: {\it  Duality for regular systems of
weights: a pr\'{e}cis.} In: {\it  Topological Field Theory, Primitive Forms and
Related Topics} (Kashiwara, M., Matsuo, A., Saito, K., Satake,I., eds.),
Progress in Math., Vol. {\bf 160}, Birkh\"auser, Boston Basel Berlin, 1998,
pp. 379--426.

\bibitem{Saito98b} Saito, K.: {\it  Duality for regular systems of
weights.} Asian J. Math. {\bf 2}, no. 4, 983--1047 (1998).

\bibitem{MSaito88} Saito, M.: {\it  Exponents and Newton polyhedra of
isolated hypersurface singularities.} Math. Ann. {\bf 281}, 411--417 (1988).

\bibitem{MSaito90} Saito, M.: {\it  Mixed Hodge modules.} Publ. Res. Inst.
Math. Sci. Kyoto Univ. {\bf 26}, 221--333 (1990).

\bibitem{MSaito93} Saito, M.: {\it  On $b$-function, spectrum and
rational singularity.}
Math. Ann. {\bf 295}, 51--74 (1993).

\bibitem{MSaito00} Saito, M.: {\it  Exponents of an irreducible plane
curve singularity.} Preprint, math.AG/0009133.

\bibitem{Scherk80} Scherk, J.: {\it  On the monodromy theorem for
isolated hypersurface singularities.} Invent. math. {\bf 58}, 289--301
(1980).

\bibitem{SchSt85} Scherk, J., Steenbrink, J.~H.~M.: {\it  On the mixed
Hodge structure on the cohomology of the Milnor fibre.} Math. Ann. {\bf
271}, 641--665 (1985).

\bibitem{Schmid73} Schmid, W.: {\it  Variation of Hodge structure: the
singularities of the
period mapping.} Invent. math. {\bf 22}, 211--319 (1973).

\bibitem{Schulze99} Schulze, M.: {\it  Computation of the monodromy of
an isolated hypersurface singularity.} Diplomarbeit, Universit\"at
Kaiserslautern, 1999 (http://www.mathematik.uni-kl.de/mschulze).

\bibitem{Schulze03} Schulze, M.: {\it  Monodromy of hypersurface
singularities.} Acta Appl. Math. {\bf 75}, 3--13 (2003).

\bibitem{SchSt00} Schulze, M., Steenbrink, J.~H.~M.: {\it  Computing
Hodge-theoretic invariants of singularities.} In: {\it  New Developments in
Singularity Theory} (Siersma, D., et al., eds.), NATO Science Series II, Vol.
{\bf 21}, Kluwer Acad. Publ., Dordrecht et al., 2001, pp. 217--233.

\bibitem{ST71} Sebastiani, M., Thom, R.: {\it  Un r\'esultat sur la monodromie.}
Invent. math. {\bf 13}, 90--96 (1971).

\bibitem{Siersma01} Siersma, D.: {\it  The vanishing topology of non
isolated singularities.} In: {\it  New Developments in
Singularity Theory} (Siersma, D., et al., eds.), NATO Science Series II, Vol.
{\bf 21}, Kluwer Acad. Publ., Dordrecht et al., 2001, pp. 447--472.

\bibitem{Steenbrink77a} Steenbrink, J.~H.~M.: {\it  Intersection form for quasi-homogeneous
singularities.} Compositio Math. {\bf 34}, 211--223 (1977).

\bibitem{Steenbrink77b} Steenbrink, J.~H.~M.: {\it  Mixed Hodge
structure on
the vanishing cohomology.} In: {\it  Real and complex singularities,} Oslo 1976 (Holm, P., ed.),
Sijthoff-Noordhoff, Alphen a/d Rijn 1977, pp. 525--563.

\bibitem{Steenbrink85} Steenbrink, J.~H.~M.: {\it  Semicontinuity of
the singularity spectrum.} Invent. Math. {\bf 79}, 557--565 (1985).

\bibitem{Steenbrink99} Steenbrink, J.~H.~M.: {\it  Spectra of $\cal
K$-unimodal isolated singularities of complete intersections.} In: {\it 
Singularity Theory,}
Proceedings of the European
Singularities Conference, Liverpool 1996 (Bruce, J.~W., Mond, D., eds.), London
Math. Soc. Lecture Note Ser., Vol. {\bf 263}, Cambridge University Press,
Cambridge, 1999, pp. 151--162.

\bibitem{Stevens03} Stevens, J.: {\it  Poincar\'e series and zeta
function for an irreducible
plane curve singularity.} Preprint, math.AG/0310215 (Bull. London Math. Soc., to appear).

\bibitem{Vaquie92} Vaqui\'e, M.: {\it  Irr\'egularit\'e des
rev\^etements cycliques.} In: {\it  Singularities} (Lille 1991) (Brasselet, J.-P.,
ed.),
London
Math. Soc. Lecture Note Ser., Vol. {\bf 201}, Cambridge
University Press, Cambridge, 1992, pp. 383--419.

\bibitem{Varchenko76} Varchenko, A.~N.: {\it  Zeta-function of
monodromy and Newton's diagram.} Invent. math. {\bf 37}, 253--262 (1976).

\bibitem{Varchenko80} Varchenko, A.~N.: {\it  Gauss-Manin connection
of isolated singular point and Bernstein polynomial.} Bull. Sci. Math. {\bf
104}, 205--223 (1980).

\bibitem{Varchenko82} Varchenko, A.~N.: {\it  Asymptotic Hodge
structure in the vanishing cohomology.} Math. U.S.S.R. Izv. {\bf 18},
469--512 (1982).

\bibitem{VKh85} Varchenko, A.~N., Khovanskii, A.~G.: {\it  Asymptotics of
integrals over vanishing cycles and the Newton polyhedron.} Soviet Math.
Dokl. {\bf 32}, 122--127 (1985).

\bibitem{Veys93} Veys, W.: {\it  Poles of Igusa's local zeta function and
monodromy.} Bull. Soc. Math. France {\bf 121}, 545--598 (1993).

\bibitem{Veys03} Veys, W.: {\it  Arc spaces, motivic integration and
stringy invariants.} Preprint, math.AG/0401374.

\bibitem{Veys04} Veys, W.: {\it  Stringy invariants of normal surfaces.}
J. Algebraic Geom. {\bf 13}, 115--141 (2004).

\bibitem{Wall83} Wall, C.~T.~C.: {\it  Classification of unimodal isolated
singularities of complete intersections.} Proc. Symp. Pure Math. Vol.  {\bf 40},
Part 2, 625--640 (1983).

\bibitem{XY89} Xu, Y., Yau, St.~S.-T.: {\it  Classification of topological
types of isolated
quasi-homogeneous two dimensional hypersurface singularities.} Manuscripta
math. {\bf 64}, 445--469
(1989).

\bibitem{Yano83} Yano, T.: {\it  $b$-functions and exponents of
hypersurface isolated singularities.} Proc. Symp. Pure Math. Vol.  {\bf 40},
Part 2, 641--652 (1983).

\bibitem{Zariski32} Zariski, O.: {\it  On the tologogy of algebroid
singularities.} Amer. J.
Math. {\bf 54}, 453--465 (1932).



\end{thebibliography}
\end{document}